\documentclass[12pt]{article}
\usepackage{amsmath,amsfonts,amssymb}
\usepackage{graphicx}

\newcommand{\fa}{\mathfrak a} 
\newcommand{\fg}{\mathfrak g}
\newcommand{\fk}{\mathfrak k}
\newcommand{\fp}{\mathfrak p}
\newcommand{\fork}{\pitchfork}
\newcommand{\acts}{\curvearrowright}
\newcommand{\into}{\hookrightarrow}
\newcommand{\bs}{\backslash}
\newcommand{\const}{\mathrm{const}}
\newcommand{\Oh}{\mathcal O}
\newcommand{\C}{\mathbb C}
\newcommand{\HH}{\mathbb H} 
\newcommand{\PP}{\mathbb P} 
\newcommand{\R}{\mathbb R}
\newcommand{\Z}{\mathbb Z}
\newcommand{\ga}{\gamma}
\newcommand{\Ga}{\Gamma}
\newcommand{\ka}{\kappa}
\newcommand{\la}{\lambda}
\newcommand{\fie}{\varphi}
\newcommand{\sM}{\mathcal M}
\newcommand{\sT}{\mathcal T}
\newcommand{\fq}{\mathfrak q}
\DeclareMathOperator{\ad}{ad}
\DeclareMathOperator{\Ad}{Ad}
\DeclareMathOperator{\Aut}{Aut}
\DeclareMathOperator{\End}{End}
\DeclareMathOperator{\Hom}{Hom}
\DeclareMathOperator{\diag}{diag}
\DeclareMathOperator{\rank}{rank}
\DeclareMathOperator{\GL}{GL} 
\DeclareMathOperator{\PSL}{PSL}
\DeclareMathOperator{\PSO}{PSO}
\DeclareMathOperator{\SL}{SL}
\DeclareMathOperator{\Spg}{Sp} 
\DeclareMathOperator{\SO}{SO}
\DeclareMathOperator{\Or}{O} 
\DeclareMathOperator{\SU}{SU}
\DeclareMathOperator{\U}{U}

\newtheorem{thm}{Theorem}
\newtheorem{conj}[thm]{Conjecture}
\newtheorem{dfn}{Definition} 
\newtheorem{exa}{Example} 
\newtheorem{exas}{Examples} 
\newtheorem{lem}{Lemma} 
\newtheorem{prob}[thm]{Problem}
\newtheorem{probA}{Problem} 
\newtheorem{probAprime}{Problem} 
\newtheorem{probC}{Problem} 

\title{On discontinuous group actions on \\ non-Riemannian homogeneous
spaces
\thanks{Translation from Sugaku, {\bf57} (2005) 267--281, to appear in
AMS translation journal ``Sugaku Expositions''}}
\author{Toshiyuki Kobayashi}
\date{}

\begin{document}
\maketitle

\noindent
\textit{Mathematics Subject Classifications} (2000):
Primary: 22F30; 

\noindent
Secondary:
  22E40,  
  53C50,  
  57S30 

\medskip
\noindent
\textit{Key words}: discontinuous group, properly discontinuous action,
discrete subgroup, proper action, Clifford--Klein form,
rigidity, deformation

\begin{abstract}
This article gives an up-to-date account of the
theory of discrete group actions on non-Riemannian
homogeneous spaces.

As an introduction of the motifs of this article,
we begin by reviewing the current knowledge of
possible global forms of pseudo-Riemannian manifolds
with constant curvatures, and
discuss what kind of problems we propose to pursue.

For pseudo-Riemannian manifolds,  
isometric actions of discrete groups are not always
properly discontinuous. 
The fundamental problem is to understand when  discrete subgroups
of Lie groups $G$ act properly discontinuously
on homogeneous spaces $G/H$ for non-compact $H$.
For this, we introduce the concepts from a
  group-theoretic perspective,
  including the `discontinuous dual' of $G/H$
  that recovers $H$ in a sense.

We then summarize recent results giving criteria
for the existence of properly discontinuous subgroups,
  and the known results and conjectures on
  the existence of cocompact ones.
The final section discusses the deformation theory
  and in particular rigidity results for cocompact
  properly discontinuous groups for pseudo-Riemannian symmetric spaces.
\end{abstract}

\section{Introduction: The problem of space forms} 

\subsection{Pseudo-Riemannian manifolds of constant curvature}

The local to global study of geometries was a major trend of 20th
century geometry, with remarkable developments achieved particularly in
Riemannian geometry. In contrast, in areas such as Lorentz geometry,
familiar to us as the space-time of relativity theory, and more
generally in pseudo-Riemannian%
\footnote{That is, having a pseudometric 2-form that is not necessarily
positive definite: whereas a Riemannian metric is given by a positive
definite 2-form at every point of a manifold, a {\em pseudo-Riemannian
metric} is the generalization obtained by replacing the positive
definite condition on the form by nondegenerate. The case of {\em
Lorentz manifolds} corresponds to the nondegenerate 2-form having
signature $(n-1,1)$.
}
geo\-metries, as well as in various other kinds of geometry (symplectic,
complex geometry,\dots), surprising little is known about global
properties of the geometry even if we impose a locally homogeneous
structure.

On the other hand, in the representation theory of Lie groups and in the
area of global analysis that applies it (noncommutative harmonic
analysis), the great trends of development throughout the 20th century
include the generalizations
\[
\begin{tabular}{rcl}
compact & $\mapsto$& noncompact \\
Riemannian & $\mapsto$& pseudo-Riemannian manifolds \\
finite & $\mapsto$& infinite dimensional representations
\end{tabular}
\]
within their terms of study, and together with these, the appearance of
ground-breaking new research methods; these moreover deepened relations
with various areas of mathematics, such as PDEs, functional analysis,
differential and algebraic geometry.

Against this background, from around the mid 1980s, I began to envisage
the possibility of developing the theory of discontinuous groups also
in the world of pseudo-Riemannian manifolds. Soon afterwards, I
succeeded in proving a necessary and sufficient condition for the
Calabi--Markus phenomenon to occur, and, stimulated by this, launched a
systematic study of the general theory of discontinuous groups of
homogeneous spaces that have good geometric structures, but are not
necessarily Riemannian manifolds; for example, semisimple symmetric
spaces or adjoint orbit spaces.

Whereas the theory of discontinuous groups of Riemannian symmetric
spaces had been the center of wide and deep developments for more than a
hundred years, at the time in the 1980s, practically no other
researchers were interested in the theory of discontinuous groups for
non-Riemannian homogeneous spaces; although the starting point was
solitary, whatever I did was a new development. After the publication of
the series of papers containing foundational results
\cite{17}--\cite{19}, \cite{28}, from the early 1990s, specialists in
other areas from France and the United States such as Benoist, Labourie,
Zimmer, Lipsman, Witte and so on eventually started to join in the study
of this problem. After this, research methods developed rapidly over the
following 10 years or so, and the ideas concerning discontinuous groups
of non-Riemannian homogeneous spaces have come to relate to many areas
of mathematics, including not only Lie group theory and discrete group
theory, but also differential geometry, algebra, ergodic theory,
mechanical systems, unitary representation theory and so on (\cite{4},
\cite{5}, \cite{9}, \cite{15}, \cite{17}, \cite{28}, \cite{34},
\cite{35}, \cite{40}, \cite{43}, \cite{46}, \cite{56}, \dots), so much
so that already no single mathematician can hope to cover them all.

For example, the recent researches of Margulis, Oh, Shalom and myself
(\cite{23}, \cite{35}, \cite{40}, \cite{46}) can be viewed as another
instance of the same trend, where the fundamental question of
understanding the distinction between a discrete group action and a
discontinuous group for a non-Riemannian homogeneous space begins to tie
in to the at first sight unrelated subject of the restriction of a
unitary representation to a noncompact subgroup.

As an introduction to this article, we review and put in order the
simplest examples (in some sense) of spaces ``of constant curvature'',
and discuss what kind of problems we propose to pursue, and what is the
current state of knowledge concerning the ``possible global forms'' of
these spaces. To set these things up precisely we recall the following
definitions.

\begin{dfn}\rm A pseudo-Riemannian manifold of constant sectional
curvature is a {\em space form}. \end{dfn}

For example, for signature $(n,0)$ (Riemannian manifolds), the sphere
$S^n$ is a space form of positive curvature, and hyperbolic space a
space form of negative curvature. For signature $(n-1,1)$ (Lorentz
manifolds), {\em de~Sitter space} is a space form of positive
curvature,%
\footnote{In Calabi and Markus \cite{8}, in connection with the use of
4-dimensional Lorentz manifold as the space-time continuum of relativity
theory, the Lorentz space form of positive curvature (that is, de~Sitter
space) is called the {\em relativistic spherical space form.}
}
{\em Minkowski space} a space form of zero curvature, and {\em anti-de~Sitter
space} a space form of negative curvature.

Here since we are interested in global properties, when we say space
form, we assume that the geometry is geodesically complete. The main
topic we consider in this section is the following question:
\begin{description}
\item[Local assumption:] among pseudo-Riemannian space forms%
\footnote{Multiplying a pseudo-Riemannian metric by $-1$ changes its
signature from $(p,q)$ to $(q,p)$, and its curvature from $\ka$ to
$-\ka$.}
of signature $(p,q)$ and curvature $\ka$,
\item[Global conclusion:] do there exist any compact examples?
\item[ ] And if so, what types of group can appear as their fundamental
groups?
\end{description}

\subsection{The two dimensional case} 
The sphere $S^2$, the torus $T^2$, and the closed Riemann surface $M_g$
of genus $g\ge2$ can be given Riemannian metrics to make them
respectively into space forms of positive, zero and negative curvature.
In other words, in two dimensional Riemannian geometry, there exists a
space form of any curvature $\ka$. The same holds in general dimensions.

However, in the case of Lorentz signature $(1,1)$, there
do not exist any compact space form with $\ka\ne0$. In fact, the sphere
$S^2$ and the Riemann surface $M_g$ with $g\ge2$ do not even admit a
Lorentz metric.%
\footnote{Any paracompact manifold admits a Riemannian structure, but the
analogous result does not hold for pseudo-Riemannian structure.
}
And if $T^2$ can be given a Lorentz metric of constant curvature $\ka$
then $\ka=0$ by the Gauss--Bonnet theorem.

\subsection{The case of positive curvature} 
Among Riemannian manifolds, the sphere $S^n$ is the typical model for
a space form of positive curvature. Conversely, the only complete space
forms with this property are $S^n$, or at most $S^n$ divided by a
suitable finite group.%
\footnote{See Wolf \cite{50} for details on what kind of finite groups
one can divide by.}
We recall two classical theorems generalizing the fact that ``a space
form of positive curvature has finite fundamental group''.

In one direction, we leave unchanged the positive definite property of
the metric form (that is, a Riemannian manifold), and perturb the curvature
(or the metric itself).

\begin{thm}[Myers 1941 \cite{38}] Suppose that the Ricci curvature of a
complete Riemannian manifold has a positive lower bound. Then the
fundamental group is finite, and the manifold is compact. \end{thm}

In the other direction, we now leave unchanged the positive definite
property of the curvature, but vary the positive definite assumption on
the metric form (that is, the condition for a Riemannian manifold).

\begin{thm}[Calabi and Markus 1962 \cite{8}] In Lorentz geometry of
dimension $\ge3$, a space form of positive curvature has finite
fundamental group, and is noncompact. \end{thm}

Theorem~2 can be formulated more generally, with {\em pseudo-Riemannian
manifold of general signature} in place of Lorentz manifold, and {\em
locally homo\-geneous space} in place of constant sectional curvature.
We can formalize this as the problem of discontinuous groups for
homogeneous spaces. Here we say that a discrete subgroup $\Ga$ of $G$ is
a discontinuous group for the homogeneous space $G/H$ to mean that the
left action of $\Ga$ on $G/H$ is properly discontinuous and free (for
more details, see Section~2). The following result is formulated
so as to contain Theorem~2 as a special case.

\begin{thm}[Criterion for the Calabi--Markus phenomenon 1989 \cite{17}]
Let $G\supset H$ be a pair of reductive Lie groups; then the homogeneous
space $G/H$ admits a discontinuous group of infinite order if and only
if\/ $\rank_\R G>\rank_\R H$. \end{thm}

In the opposite directions, testing this theorem on various examples of
homogenous spaces leads one to believe that the following conjecture may
hold:

\begin{conj}[see \cite{27}] Assume that $p\ge q>0$ and $p+q\ge3$.
Suppose that the sectional curvature of a complete pseudo-Riemannian
manifold of signature $(p,q)$ has a positive lower bound.%
\footnote{Here we are using sectional curvature. With Ricci curvature
only, the assumption is too weak.}
Then the fundamental group is finite and the manifold is noncompact.
\end{conj}

\subsection{The case of curvature zero} 
Among Riemannian manifolds, the $n$-dimensional torus $T^n$ is
the typical example of a space form of curvature 0. Its fundamental
group $\Z^n$ is an Abelian group. More generally, the following theorem
says that the fundamental group of a space form of curvature~0 is also
close to an Abelian group:

\begin{thm}[Bieberbach 1911] The fundamental group of a complete
Riemannian manifold of constant sectional curvature~$0$ contains an
Abelian subgroup of finite index. \end{thm}

Whether the analogous result holds for a pseudo-Riemannian manifold is
still unknown:

\begin{conj}[Auslander Conjecture -- Special case] The fundamental group
of a compact pseudo-Riemannian manifold of constant sectional
curvature~$0$ contains a solvable subgroup of finite index. \end{conj}

Conjecture~6 is true for a Lorentz manifold (Goldman and
Kamishima 1984 \cite{12}, Tomanov). More generally, one can extend the
Bieberbach theorem under the assumption of an affine manifold. This is
the original Auslander Conjecture. One can also envisage a stronger form:

\begin{prob}[Milnor 1977 \cite{37}] \rm Is it true that any
discontinuous group (see Section~2) for affine space
$\R^n=(\GL(n,\R)\ltimes\R^n)/\GL(n,\R)$ contains a solvable subgroup of
finite index? \end{prob}

In 1983 Margulis gave a counterexample to this conjecture of Milnor in
dimension $n=3$ (Theorem~11). On the other hand, the following
proposition is a continuous analog of Milnor's problem:
\begin{quote}
A connected subgroup of the affine group that acts properly on $\R^n$
(see Section~2) is amenable; that is, it is a compact extension of a
solvable group.
\end{quote}
This  statement (also in a more general form) is known to hold (1993 \cite{20}, Lipsman 1995
\cite{34}). Although the original Auslander Conjecture remains open,
Abels, Margulis and Soifer (1997 \cite{1}, \cite{2}, \cite{36}) have
announced that it holds in dimension $\le6$. Also, on a related topic,
the Lipsman conjecture (1995, \cite{34}) on discontinuous groups of
nilmanifolds and their proper actions is known. Definitive results on
the Lipsman conjecture have been obtained in very recent work of
Nasrin, Yoshino, Baklouti and Khlif. Namely, it is true if the
nilmanifold has dimension~$\le4$ \cite{51}, and there is a
counterexamples in dimension $\ge5$ \cite{52}. Moreover, it is
true for Lie groups that are at most 3-step nilpotent (\cite{3},
\cite{39}, \cite{54}), and there is a counterexample for 4-step
nilpotent or more \cite{52}.

\subsection{The case of negative curvature} 
In the case of Riemannian manifolds, negatively curved compact space
forms (hyperbolic manifolds) exist. This is equivalent to the fact that
the Lorentz groups $\Or(n,1)$ admit uniform lattices.%
\footnote{For arithmetic uniform lattices, this is general theory, due
to Borel and Harish-Chandra (1962), Mostow and Tamagawa (1962) and Borel
(1963); for nonarithmetic uniform lattices (in hyperbolic spaces)
examples were constructed by Makarov (1966), Vinberg, Gromov and
Pyatetski{\u i}-Shapiro (1981).
}
However, for general pseudo-Riemannian manifolds, the fundamental
problem of knowing for which signature $(p,q)$ there exist {\em compact
space forms (of negative curvature)} is still not completely settled as
things stand. The following conjecture addresses this question; it is a
specialization to the case of the homogeneous spaces
$\Or(p,q+1)/\Or(p,q)$ of Conjecture~17, the existence problem for
uniform lattices for a reductive homogeneous space.

\begin{conj}[Conjecture on Space Forms, 1996]
The necessary and sufficient condition for the existence of a compact
pseudo-Riemannian manifold of signature $(p,q)$ with constant negative
sectional curvature is that $(p,q)$ is in the following list:
\[
\begin{array}{|c||c|c|c|c|c|}
\hline
q&N&0&1&3&7 \\
\hline
p&0&N&2N&4N&8 \\
\hline
\end{array}
\]
where $N=1,2,3,\dots$.
\end{conj}

The sufficiency of the condition is proved. As already discussed, for
$q=0$, these are hyperbolic (Riemannian) manifolds; for $q=1,3$,
examples were discovered in 1981 by Kulkarni \cite{30}. In the case
$q=7$, examples were discovered from the 1990s by applying the existence
theorem of uniform lattices for general homogeneous spaces (Theorem~15),
see \cite{21}.

Whether the condition of Conjecture~8 is necessary is still not settled,
although it is proved to hold in many cases, such as $q=1$, or $p\le q$, or
$pq$ odd. The final ``odd condition'' on $p,q$ was extended to general
reductive homogeneous spaces (the most typical pseudo-Riemannian
homogeneous spaces) by Ono and myself \cite{28}, by a method
generalizing Hirzebruch proportionality for characteristic classes. (For
more details on these topics, we refer to the references given in
\cite{27} and \cite{31}. See also \cite{KY}.)

\section{Discontinuous actions and Clifford--Klein forms} 
Even though the problem of space forms treats extremely special spaces,
as discussed in the preceding section, many problems remain open.
However, even restricting to these cases, there are instances when,
rather than studying individual isolated examples, we obtain a clearer
perspective from the general viewpoint of {\em discontinuous group
actions on (non-Riemannian) homogeneous spaces.}%
\footnote{From this point of view, any pseudo-Riemannian space forms of
signature $(p,q)$ with $q\ne1$ is a Clifford--Klein form of the rank~1
semisimple symmetric space $\Or(p,q+1)/\Or(p,q)$.
}
In this direction, this section explains basic notions and concrete
examples, while emphasizing the distinction between ``discontinuous
groups acting on homogeneous spaces'' and ``discrete subgroups''. In the
remainder of the article, keeping at the back of our minds the point of
view on the problem of space forms explained above, we want to discuss
how far the world of discontinuous groups can be extended in the general
framework, avoiding as far as possible the technical terms of the theory
of Lie group.

First, in the case of Riemannian manifolds, subgroups consisting of
isometries satisfy
\[
\hbox{discrete group}
\enspace \iff \enspace
\hbox{discontinuous group.}
\]
However, for pseudo-Riemannian manifolds and for subgroups consisting of
isometries, these conditions are not equivalent:
\[
\hbox{discrete group}
\enspace \stackrel
{\displaystyle \Longrightarrow\kern-6.5mm\times\kern2.2mm}
{\Longleftarrow}
\enspace
\hbox{discontinuous group,}
\]
and the quotient space by the action of a discrete group is not
necessarily Hausdorff. For example, the orbit of a discrete group is not
necessarily a closed set; this corresponds to the fact that in the
topology of the quotient space, a single point is not a closed point,
and in particular the quotient topology is not Hausdorff. There are thus
cases when the quotient is non-Hausdorff for local reasons. However,
Hausdorff is a global property of a topological space, and there are
also more curious counterexamples. The following example is one such.
Here the quotient topology is non-Hausdorff for a global reason, because
accumulation points do not exist.%
\footnote{This example can be reformulated in group theoretic terms as a
homogeneous space of $\SL(2,\R)$.
}

\begin{exa}\rm Make the discrete group $\Z$ act on
$X=\R^2\setminus\{(0,0)\}$ by the map $\Z\times X\to X$ given by
$(n,(x,y))\mapsto(2^nx,2^{-n}y)$. Thus the $\Z$-orbits are contained in
the hyperbolas $xy=\const.$ or in the $x$- and $y$-axes;
\begin{figure}[h]
\hspace*{11.5em}
\includegraphics{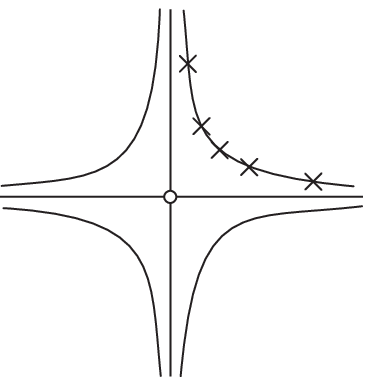}
\caption{$\Z$-orbits}
\end{figure}
in Figure~1 the crosses $\times$ represent a $\Z$-orbit contained in the
first quadrant. This $\Z$-action does not have any accumulation points
in $X$, but the quotient space $\Z\bs X$ is non-Hausdorff. In fact, by
considering the fiber bundle
$\Z\bs\R\to\Z\bs X\to\R\bs X$, one sees that the quotient space $\Z\bs
X$ is homeomorphic to an $S^1$ fiber bundle over the base space
illustrated in Figure~2,
\begin{figure}[h]
\includegraphics{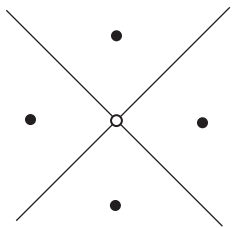}
\caption{$\R\bs X \simeq$ four half-lines and four points \newline (a
complete system of representatives)}
\end{figure}
consisting of four half-lines and four points, given a non-Hausdorff
topology. \end{exa}

How to understand this kind of example in group theoretic terms is the
main topic of Sections~2--3. As a preparation, we introduce some pieces
of basic terminology. The set-up we consider is a topological group
$\Ga$ with a continuous action on a space $X$; we write $\Ga\acts X$ for
this action. For $S\subset X$ a subset, we define the subset
$\Ga_S\subset\Ga$ as follows:
\[
\Ga_S := \bigl\{ \ga\in\Ga \bigm| \ga\cdot S\cap S \ne \emptyset \bigr\}.
\]

\begin{dfn}\rm
\begin{enumerate}
\renewcommand{\labelenumi}{(\arabic{enumi})}
\item The action $\Ga\acts X$ is {\em properly discontinuous} if $\Ga_S$
is finite for any compact subset $S\subset X$.
\item The action $\Ga\acts X$ is {\em proper} if $\Ga_S$ is compact for
any compact subset $S\subset X$.
\item The action $\Ga\acts X$ is {\em free} if the stabilizer subgroup
$\Ga_{\{p\}}$ of any point $p\in X$ consists of the identity element
only.
\end{enumerate}
\end{dfn}

The action of a noncompact group is not necessarily well behaved. The
notion of a ``proper action'' (Palais 1961 \cite{44}) abstracts out the
``good property'' of a compact group action. Putting together as
above the notion of properly discontinuous gives rise to the following
equation:
\[
\hbox{$\Ga$ acts properly discontinuously} \iff\enspace
\left\{ \begin{array}{l}
\hbox{$\Ga$ acts properly,} \\
\hbox{{\em and\/} $\Ga$ is a discrete group.}
\end{array} \right.
\]
Thus, in order to determine whether the action of a discrete group is
properly discontinuous, it is enough to determine whether the action is
proper. This last point has wide applications.

Now, for a group $\Ga$ acting on a set $X$, we write $\Ga\bs X$ for the
set of equivalence classes of the equivalence relation given by
\[
x\sim x' \iff \hbox{$x'=\ga\cdot x$ for some $\ga\in\Ga$.}
\]
We can view $\Ga\bs X$ as the set of $\Ga$-orbits in $X$, so we call it
the {\em orbit space} (or $\Ga$-orbit space). The reason for considering
properly discontinuous actions of $\Ga$ is the following well-known
result:

\begin{lem} Let $X$ be a manifold (respectively $C^\infty$ manifold,
pseudo-Riemannian manifold, complex manifold, etc.), and suppose that a
discrete group $\Ga$ acts on $X$ continuously (respectively smoothly,
isometrically, bi\-holo\-morphically, etc.). If the action of $\Ga$ is
properly discontinuous and free then the quotient topology on $\Ga\bs X$
is Hausdorff, and the quotient $\Ga\bs X$ can be given a unique manifold
structure for which the quotient map $X\to\Ga\bs X$ is a local
homeomorphism (respectively local diffeomorphism, local pseudo-isometry,
locally biholomorphic). \end{lem}

In what follows, $X=G/H$ is a homogeneous space for a Lie group $G$, and
$\Ga$ is a discrete subgroup of $G$, so we have the following triple of
groups:
\[
\Ga\subset G\supset H.
\]

\begin{dfn}\rm We say that $\Ga$ is a {\em discontinuous group} for the
homogeneous space $G/H$ if $\Ga$ acts properly discontinuously and
freely on $G/H$. Here ``properly discontinuously'' is the most important
condition; papers in the literature sometimes omit the condition that the
action is free in the definition of a discontinuous group. If $\Ga$ is a
discontinuous group for $G/H$, the manifold obtained as the double coset
space $\Ga\bs G/H$ is a {\em Clifford--Klein form} of $G/H$. If in
addition $\Ga\bs G/H$ is compact, we say that $\Ga$ is a {\em uniform
lattice} of the homogeneous space $G/H$. \end{dfn}

\begin{exas}\rm
\begin{enumerate}
\renewcommand{\labelenumi}{(\arabic{enumi})}
\item Let $(G,\Ga,H)=(\R^n,\Z^n,\{0\})$; then the Clifford--Klein form
$\Ga\bs G/H$ is diffeomorphic to the $n$-dimensional torus $T^n$, and so
is a compact manifold.

\item Nilmanifold: set
\[
G:=\left.\left\{\begin{pmatrix} 1&a&b\\ 0&1&c\\ 0&0&1\end{pmatrix}
\right|a,b,c,\in\R\right\},
\quad \Ga := G\cap\GL(3,\Z),
\quad H = \{e\}.
\]
The Clifford--Klein form $\Ga\bs G/H$ is a 3-dimensional compact
manifold (the {\em Iwasawa manifold}).

\item The modular group: set $(G,\Ga,H)=(\SL(2,\R),\SL(2,\Z),\{e\})$.
The Clifford--Klein form $\Ga\bs G/H$ is noncompact, but has finite
volume (with respect to a naturally defined measure). Moreover $\Ga\bs
G/H$ is homeo\-morphic to the knot complement
$\R^3\setminus\hbox{(trefoil knot)}$.

\item Closed Riemann surface: the closed Riemann surface $M_g$ of genus
$g\ge2$ can be realized as a Clifford--Klein form $\Ga\bs G/H$ of the
Poincar\'e upper half-space $G/H=\SL(2,\R)/\SO(2)$; here
$\Ga\simeq\pi_1(M_g)$.

\item The Calabi--Markus phenomenon: let $G=\SL(2,\R)$, and let
$H\subset G$ be any noncompact closed subgroup. Then the only
discontinuous groups of $G/H$ are finite groups. In particular, if $H$
and $G/H$ are both noncompact then there does not exist any uniform
lattice for $G/H$.

\item Compact Lorentz space form: let $(G,H)=(\SO(2n,2),\SO(2n,1))$ and
let $\Ga\subset\U(n,1)$ be a uniform lattice without torsion elements.
If we view $\Ga$ as a discrete subgroup of $G$ by the embedding
$\Ga\subset\U(n,1)\subset G$ then $\Ga$ is also a uniform lattice of the
homogeneous space $G/H$. However, $\Ga$ cannot be a uniform lattice of
$G$.
\end{enumerate}
\end{exas}

The above examples (5) and (6) illustrate the following important
warning for noncompact subgroups:
\begin{quote}
for a noncompact subgroup $H$, a uniform lattice for $G$ is {\em not}
the same thing as a uniform lattice for the homogeneous space $G/H$.
\end{quote}

\section{Criterion for an action to be discontinuous} 
This section discusses the following problem:

\begin{probA} \rm Find effective methods of determining whether a discrete
subgroup $\Ga$ acts properly discontinuously on a homogeneous space
$G/H$. \end{probA}

The definition of a properly discontinuous action on a topological space
was easy enough to formulate. However, in general, given a discrete
subgroup $\Ga$ of a Lie group $G$, actually determining whether or not
the action of $\Ga$ on a homogeneous space $G/H$ is properly
discontinuous is not at all easy. One aims for ``criteria'' in
Problem~A that are so concrete and powerful that we can, for example,
obtain the following various theorems as sample applications.

\begin{thm}{\bf(Pseudo-Riemannian manifold space form of signature
$(p,q)$ with negative curvature \cite{8}, \cite{48}, \cite{30},
\cite{17})}\enspace
The homo\-geneous space $\Or(p,q+1)/\Or(p,q)$ has only finite groups as
discontinuous groups $\iff p\le q$. \end{thm}

\begin{thm}[Solvable manifolds 1993 \cite{20}, \cite{34}] Any
homogeneous space for a solvable Lie group has a Clifford--Klein form
with fundamental group of infinite order. \end{thm}

\begin{thm}[Affinely flat manifold, Margulis 1983 \cite{2}] Affine
$3$-space $(\GL(3,\R)\ltimes\R^3)/\GL(3,\R)\simeq\R^3$ admits a free
non-Abelian group as a discontinuous group. \end{thm}

\begin{thm}[Pseudo-Riemannian homogeneous space, Benoist 1996 \cite{4}]
$\SL(3,\R)/\SL(2,\R)$ does not admit a free non-Abelian group as a
discontinuous group. \end{thm}

Now for a non-Riemannian homogeneous space $G/H$, the usual approach to
its discontinuous groups was to restrict the study to extremely special
situations (for example, the rank~1 symmetric spaces of the type treated
in Section~1) and to make clever use of the special properties enjoyed
by the individual homogeneous spaces $G/H$; compare \cite{8}, \cite{30},
\cite{48}, \cite{49},\dots. This method requires huge calculations, even
if $G/H$ is a rank~1 symmetric space (as in \cite{30}). Instead of this,
to deal with discontinuous groups for more general pseudo-Riemannian
homogeneous spaces (for noncompact $H$), I introduced the following idea
in \cite{17}, \cite{22}.

\begin{enumerate}
\item {\em Forget} that $H$ is a group and that the homogeneous space $G/H$ is
a manifold.
\item {\em Forget} that $\Ga$ is discrete and that it is a group.
\end{enumerate}

Having thus thrown away all of the (at first sight) most important
information, we are left with the following possibilities:

\begin{enumerate}\setcounter{enumi}{2}
\item View $\Ga$ and $H$ {\em on an equal footing}, simply as subsets of
$G$.
\item Control the discontinuous property of the action $\Ga\acts G/H$
using the representation theory of $G$.
\end{enumerate}

In order to implement this idea, we introduce the following two
relations $\fork$ and $\sim$ on subsets $H$ and%
\footnote{In what follows, we often use $L$ as an alternative notation
to $\Ga$.}
$L$ of a locally compact group $G$.

\begin{dfn}[see \cite{22}]\rm \begin{enumerate}
\renewcommand{\labelenumi}{(\arabic{enumi})}
\item We say that the pair $(L,H)$ is {\em proper} in $G$ and write
$L\fork H$ (in $G$) if and only if for any compact subset $S$ of $G$ the
intersection $L\cap SHS$ is relatively compact.%
\footnote{In differential geometry, $\fork$ is often used to mean that
two submanifolds intersect transversally; here we use the notation in a
completely different meaning.
}

\item We write $L\sim H$ (in $G$) if and only if there exists a compact
subset $S$ of $G$ such that $L\subset SHS$ and $H\subset SLS$.
\end{enumerate}
\end{dfn}
For an Abelian group $G$, the relations $\fork$ and $\sim$ are
remarkably simple.

\begin{exa}\rm Let $H$ and $L$ be subspaces of the Abelian group
$G:=\R^n$.
\begin{enumerate}
\renewcommand{\labelenumi}{(\arabic{enumi})}
\item $H\fork L \hbox{ in $G$} \iff H\cap L=\{0\}$.
\item $H\sim L \hbox{ in $G$} \iff H=L$.
\end{enumerate}
\end{exa}

Now if $L$ and $H$ are closed subgroups of $G$ then
\[
\hbox{$L\fork H$ in $G$} \iff \hbox{$L$ acts properly on $G/H$.}
\]
Thus we can view $\fork$ as a notion generalizing properness of a group
action. In other words, to understand whether an action is proper, or is
properly discontinuous, it is enough to understand the relation $\fork$.
Moreover
\[
L\fork H \iff H\fork L
\]
This reflects a kind of symmetry between the action of $L$ on $G/H$ and
of $H$ on $G/L$. We also have
\[
\hbox{if $H\sim H'$ then}\quad L\fork H \iff L\fork H'.
\]
Thus the use of $\sim$ provides economies in considering $\fork$. We
define the {\em discontinuous dual} $H^\fork$ of a subset $H$ of $G$ as
follows:
\[
H^\fork := \bigl\{ L \bigm| \hbox{$L$ is a subset of $G$ satisfying
$L\fork H$} \bigr\}.
\]

First, we have the following theorem, which I proved in 1996
\cite{22} for $G$ a reductive Lie group; whether it holds in general was
one of the unsolved problems discussed in \cite{27}, and was settled
positively by Yoshino in 2004 \cite{53}.

\begin{thm}[Duality theorem] A subset $H$ of a Lie group $G$ can be
reconstructed from its discontinuous dual $H^\fork$ up to the
equivalence $\sim$. \end{thm}

Our original aim was to determine by explicit methods whether the action
of a discrete group is properly discontinuous. Problem~A can be
formalized again in the following more general form.

\begin{probAprime} \rm Let $H$ and $L$ be subsets of a group $G$. Find
criteria to determine whether $H\fork L$. \end{probAprime}

Let $G$ be a reductive linear Lie group (for example, $\GL(n,\R)$ or
$\Or(p,q)$, etc.). Write $\fg=\fk+\fp$ for a Cartan decomposition of the
Lie algebra $\fg$ of $G$, and choose a maximal Abelian subspace $\fa$ of
$\fp$. Write
\[
d(G) := \dim\fp,\quad \rank_\R G :=\dim a.
\]
Also, using a Cartan decomposition $G=K\exp(\fa)K$, we define the {\em
Cartan projection} $\nu\colon G\to\fa$, which is determined up to the
action of the Weyl group.

For example, for $G=\GL(n,\R)$, we have $d(G)=\binom{n+1}2$ and
$\rank_\R G=n$. For a square matrix $g\in G$, the product $^tgg$ is a
positive definite symmetric matrix, and we write out its eigenvalues in
order, from the largest down:
\[
\la_1\ge\cdots\ge\la_n \ (>0).
\]
Then the Cartan projection $\nu\colon G\to\R^n$ is given by the formula:
$g\mapsto \frac12(\log\la_1,\dots,\log \la_n)$.

After these preparations, the answer to Problem~A (or Problem~A$'$) is
as follows:

\begin{thm}[Criterion for a properly discontinuous action \cite{17},
\cite{22}] Let $H$ and $L$ be subsets of a reductive linear Lie group
$G$. Then
\begin{enumerate}
\renewcommand{\labelenumi}{$(\arabic{enumi})$}
\item $L\sim H$ in $G$ $\iff$ $\nu(L)\sim\nu(H)$ in $\fa$.
\item $L\fork H$ in $G$ $\iff$ $\nu(L)\fork\nu(H)$ in $\fa$.
\end{enumerate}
\end{thm}

For the Abelian group $\fa\simeq\R^n$ the relations $\fork$ and $\sim$
have a very simple meaning. Thus Theorem~14 is useful as a criterion.

I solved Problem~A$'$ in 1989 \cite{17} for a triple of reductive Lie
groups $(G,H,L)$, and subsequently generalized the result in \cite{22}
in the above form, without even assuming that $L$ and $H$ have group
structures; Benoist \cite{4} proved a similar generalization
independently of \cite{22}. The above results Theorems~3, 9 and~12 can
be obtained as corollaries of Theorem~14. Also, using a generalization
in the style of Theorem~14 made it possible to study the extent to which
the proper discontinuity of an action is preserved on deforming a
discrete group (for example, Goldman's conjecture 1985 \cite{11}); see
\cite{24}, \cite{29}, \cite{41}, \cite{45}; we return to this topic in
Section~5. In addition to this, recent work (from 2001 onwards) due to
Iozzi, Oh, Witte and others \cite{14}, \cite{27}, \cite{42} contain
results arising from applications of Theorem~14 to individual
homogeneous spaces.

The implications $\Leftarrow$ of Theorem~14, (1) and $\Rightarrow$ of
Theorem~14, (2) are obvious. And the implication $\Rightarrow$ in (1) is
related to giving uniform bounds on the errors in the eigenvalues when a
matrix is perturbed.%
\footnote{Various inequalities concerning this are known, of which a
theorem of Weyl is the prototype.
}

\section{The existence problem for compact Clifford--Klein forms}
\subsection{Existence and nonexistence theorems for uniform lattices}

In what follows, we let $G\supset H$ be a pair of reductive linear Lie
groups. $G/H$ is the typical model of a pseudo-Riemannian homogeneous
space (Riemannian if $H$ is compact). This section discusses the
following problem.

\begin{probA} \rm For which homogeneous spaces $G/H$ does a uniform lattice
exist? In other words, classify the homogeneous spaces that have compact
Clifford--Klein forms. \end{probA}

Among the various currently known results, we discuss two that have the
widest field of applications. The key points in the proof of these
results are the criterion for a proper action (Theorem~14) and the
cohomology of discrete groups.

\begin{thm}[1989 \cite{17}] If a reductive subgroup $L\subset G$
satisfies
\[
L\fork H \quad\hbox{and}\quad d(L)+d(H)=d(G),
\]
then a compact Clifford--Klein form of $G/H$ exists.
\end{thm}

\begin{thm}[1992 \cite{19}] If $L\subset G$ is a reductive subgroup and
there exists $H$ with
\[
L\sim H \quad\hbox{and}\quad d(L)>d(H)
\]
then there does not exist a compact Clifford--Klein form of $G/H$.
\end{thm}

We refer to \cite{21} for the list of pairs $(G,H)$ that satisfy the
conditions of Theorem~15%
\footnote{It follows of course from Borel's Theorem \cite{5} that this
list includes the case that $H$ is compact, and the case of the group
manifold itself (that is $G=G' \times G'$ and $H=\diag G'$).
}
and to \cite{19}, \cite{21} for the list of $(G,H)$ that satisfy those of
Theorem~16.

If $\Ga$ is a torsion-free uniform lattice of a group $L$ then, provided
the assumptions of Theorem~15 are satisfied, $\Ga\bs G/H$ is a compact
Clifford--Klein form. Conversely, even if we assume that $\Ga\bs G/H$ is
a compact Clifford--Klein form, it does not necessarily follow that
there exists a reductive Lie subgroup $L$ containing $\Ga$ satisfying
the assumptions of Theorem~15 (1998 \cite{24}, Salein 1999 \cite{45}).
However, in a slightly weaker form, the following conjecture is still
unsolved. The special case $G/H=\Or(p,q+1)/\Or(p,q)$ of Conjecture~17 is
Conjecture~6 on space forms.

\begin{conj}[1996 \cite{21}] The converse of Theorem~15 also holds.
\end{conj}

\subsection{Uniform lattices of an adjoint orbit} 
We consider semisimple orbits as examples of homogeneous spaces. If we
choose an element $X$ of the Lie algebra $\fg$ of a Lie group $G$, the
adjoint orbit $\Oh_X =\Ad(G)X$ is a submanifold of $\fg$ that we can
identify with the homogeneous space $G/G_X$, where $G_X=\{g\in
G\bigm|\Ad(g)X=X\}$ is the stabilizer of $X$.

When $G$ is a reductive Lie group and $\ad X\in\End(\fg)$ is semisimple,
we say that $\Oh_X$ is a {\em semisimple orbit}. If in addition all of
the eigenvalues of $\ad X$ are purely imaginary numbers then we say that
$\Oh_X$ is an {\em elliptic orbit}. For example, the adjoint orbits of a
compact Lie group are always elliptic orbits.

A number of important classes of homogeneous spaces arise as semisimple
orbits. For example, all the flag manifolds, Hermitian symmetric spaces,
para-Hermitian symmetric spaces, and so on can be realized as semisimple
orbits.%
\footnote{The first two cases are even elliptic orbits.
}

One can define a natural $G$-invariant symplectic structure and
pseudo-Riemannian structure on a semisimple orbit.%
\footnote{Plus, in the case of an elliptic orbit, a $G$-invariant
complex structure. A K\"ahler or pseudo-K\"ahler structure can also be
defined.
}
Also, various unitary representations appear as geometric quantizations
of these orbits, including principal series representations (more
generally, degenerate principal series representations), and discrete
series representations (more generally, Zuckerman's derived functor
modules $A_\fq(\la)$), and so on. The latter in particular corresponds
to geometric quantization of elliptic orbits.

We have the following theorem concerning the existence problem%
\footnote{Although extremely special, the case that a compact
Clifford--Klein form exists has known analytic applications, such as
Atiyah and Schmid's construction \cite{AS} of the discrete series
representations via the $L^2$ index theorem.
}
of compact Clifford--Klein forms of semisimple orbits.

\begin{thm}[see \cite{19}] The only semisimple orbits having a uniform
lattice are elliptic orbits. In particular, these have an invariant
complex structure. \end{thm}

Hermitian symmetric spaces always have a uniform lattice (Borel
\cite{6}), and are elliptic orbits. We give one example of an elliptic
orbit that is not a Hermitian symmetric space: consider the Hermitian
form of signature $(2,n)$
\[
z_1\overline{z_1} + z_2\overline{z_2} - z_3\overline{z_3} - \cdots
-z_{n+2}\overline{z_{n+2}},
\]
and write $\Oh\subset\PP^{n+1}_\C$ for the set of all complex lines such
that the restriction of the form is positive definite. Then $\Oh$ is an
open subset of complex projective space $\PP^{n+1}_\C$ (so is in
particular a complex manifold), and can be identified with an elliptic
orbit of $\U(2,n)$. Note that $\Oh$ does not have the structure of a
Hermitian symmetric space. Moreover, if $n$ is even then we can use
Theorem~15 to see that there exists a uniform lattice for $\Oh$ (for
this, one need only set $L=\Spg(1,\frac n2)$). In particular, one can use
this to construct a compact symplectic complex manifold for which the
natural form has indefinite signature (see \cite{17}, \cite{21}).

Theorem~18 was discovered by myself; its proof uses the cohomology of
discrete groups, and is based on Theorem~16 (\cite{18}, \cite{19}).
Subsequently Benoist and Labourie \cite{5} gave a different proof using
symplectic geometry.

\subsection{Uniform lattices of $\SL(n)/\SL(m)$} 
This section discusses the question of whether there exist compact
Clifford--Klein forms of the non-symmetric homogeneous space
$\SL(n)/\SL(m)$. This space is special from our point of view; from the
mid-1990s onwards, the question of the existence of its compact
Clifford--Klein forms was attacked using a variety of different methods,
with the same result being obtained using many different methods of
proof, resulting in an attractive amalgam with other areas.

The original model is the following result, obtained by applying
$L=\SU(2,1)$ to Theorem~16.

\begin{thm}[1990 \cite{18}] There do not exist any compact
Clifford--Klein forms of $\SL(3,\C)/\SL(2,\C)$. \end{thm}

The following theorem is deduced from Theorem~16 by the same principle
(replacing $\R$ by $\C$ or by the quaternions $\HH$ leads to similar
results).

\begin{thm}[1992 \cite{19}] There do not exist any compact
Clifford--Klein forms of the homogeneous spaces $\SL(n,\R)/\SL(m,\R)$ if
$\frac n3>\left[\frac{m+1}2\right]$. \end{thm}

Now for the Clifford--Klein forms of $\SL(n,\R)/\SL(m,\R)$, we can also
consider the right action of $\SL(n-m,\R)$. Taking note of this point,
Zimmer and his collaborators used machinery such as the superrigidity
theorem for cocycles and Ratner's theorem on the closure of orbits to
prove the following theorems.

\begin{thm}[Zimmer 1994 \cite{56}] There do not exist any compact
Clifford--Klein forms of $\SL(n,\R)/\SL(m,\R)$ when $n>2m$. \end{thm}

\begin{thm}[Labourie, Mozes and Zimmer 1995 \cite{32}] There do not
exist any compact Clifford--Klein forms of $\SL(n,\R)/\SL(m,\R)$ if
$n\ge2m$. \end{thm}

\begin{thm}[Corlette and Zimmer 1994 \cite{9}, \cite{10}] There do not
exist any compact Clifford--Klein forms of $\Spg(n,2)/\Spg(m,1)$ if
$n>2m$. \end{thm}

The following result was obtained as an application of the criterion for
a proper action (Theorem~14).

\begin{thm}[Benoist 1996 \cite{4}] There do not exist any compact
Clifford--Klein forms of $\SL(n,\R)/\SL(n-1,\R)$ for odd $n$. \end{thm}

Moreover, for $\SL(m)$ embedded in $\SL(n)$ by an irreducible
representation $\fie$ (not just the natural inclusion), Margulis
considered the restriction of unitary representations of $\SL(n,\R)$ to
noncompact subgroups, and used methods involving delicate estimates for
the asymptotic behavior of matrix coefficients to prove the following
theorem.

\begin{thm}[Margulis 1997 \cite{35}] For $n\ge5$, there does not exist
any compact Clifford--Klein form of $\SL(n,\R)/\fie(\SL(2,\R))$.
\end{thm}

The systematic study of Margulis' method was taken further by Oh 1998
\cite{40}. The method based on unitary representation theory was
developed futher, and Shalom proved the following theorem.

\begin{thm}[Shalom 2000 \cite{46}] For $n\ge4$, there does not exist any
compact Clifford--Klein form of $\SL(n,\R)/\SL(2,\R)$. \end{thm}

These results were obtained by taking up recent development in other
areas of mathematics, and the methods of proof extend over many
branches. However, as things stand at present, all the currently known
results support Conjecture~17 (which, applied to this case, states that
``there do not exist any compact Clifford--Klein forms of
$\SL(n,\R)/\SL(m,\R)$ for $n>m$'').

Note that, among these theorems, Theorems~21, 22, 23 and~26 are contained in
an extremely special case of Theorem~16: although the references \cite{9},
\cite{10}, \cite{31}, \cite{46}, \cite{56} mentioned above do not refer
explicitly to this, Theorem~16 or its corollary Theorem~20 is actually a
stronger result even when restricted to these special cases. On the
other hand, Benoist's Theorem~24 and Margulis' Theorem~25 are not
contained in Theorem~16. For various results concerning these explicit
kinds of homogeneous spaces, many details are contained in my lecture
notes \cite{21} in the proceedings of a European School.

\section{Rigidity and deformations of Clifford--Klein forms} 
This section discusses the following problem.

\begin{probA} \rm Is it possible to deform a uniform lattice $\Ga$ for a
homogenous space $G/H$? \end{probA}

For an irreducible Riemannian symmetric space $G/H$ of dimension~$\ge3$
with a compact subgroup $H$, and $\Ga$ a uniform lattice of $G/H$, there
do not exist any essential deformations of $\Ga$ (Theorem~27). This
result is the original model for various kinds of rigidity theorems (in
Riemannian geometry).

Now, does there exist a similar {\em rigidity result} in the case that
$H$ is noncompact (the pseudo-Riemannian case)? We can view a ``rigidity
theorem'' as an assertion of the type that the fundamental group
determines not just the topological structure, but also the geometric
structure. Now, does the ``rigidity theorem'' also hold for an
(irreducible) pseudo-Riemannian symmetric space?

In fact, for a noncompact subgroup $H$, the situation for the rigidity
theorem is quite different from the case of Riemannian symmetric spaces.
More precisely, there exist (irreducible) pseudo-Riemannian symmetric
spaces of arbitrarily high dimension, that admit uniform lattices for
which the rigidity theorem does not hold (see \cite{20}, \cite{24}). We
give such examples in Theorem~28, where we give a more precise
formulation of Problem~C. But first, we note that Problem~C includes the
following two subproblems.

\begin{probC} \rm For a discrete subgroup $\Ga\subset G$, describe the
deformations of $\Ga$ as an abstract group inside $G$. \end{probC}

\begin{probC} \rm If a discrete subgroup $\Ga\subset G$ can be deformed
inside $G$, determine the range of the deformation parameters that does
not destroy the proper discontinuity of its action on $G/H$. \end{probC}

Bearing these questions in mind, we now try to describe abstractly the
set of deformations of a discontinuous group.

Let $G$ be a Lie group and $\Ga$ a finitely generated group; we give the
set $\Hom(\Ga,G)$ of group homomorphisms from $\Ga$ to $G$ the pointwise
convergence topology. The same topology is obtained by taking
generators $\ga_1,\dots,\ga_k$ of $\Ga$, then using the injective map
\[
\Hom(\Ga,G)\into G\times\cdots\times G \quad\hbox{given by}\quad
\fie\mapsto(\fie(\ga_1),\dots,\fie(\ga_k))
\]
to give $\Hom(\Ga,G)$ the relative topology induced from the direct
product $G\times\cdots\times G$.

Let $H$ be a closed subgroup of $G$. As already explained, if $H$ is
noncompact then a discrete subgroup of $G$ does not necessarily act
properly discontinuously on $G/H$. Here rather than $\Hom(\Ga,G)$, it is
the subset $R(\Ga,G,H)$ defined below that plays the important role (see
\cite{20}):
\[
R(\Ga,G,H) := \left\{ u \in \Hom(\Ga,G) \left|
\renewcommand{\arraycolsep}{.25em}
\begin{array}l \hbox{$u$ is injective; and $u(\Ga)$ acts properly}\\
\hbox{discontinuously and freely on $G/H$}
\end{array} \right\}\right..
\]

Then for each $u\in R(\Ga,G,H)$ we obtain a Clifford--Klein form
$u(\Ga)\bs G/H$. Now the direct product group $\Aut(\Ga)\times G$ acts
naturally on $\Hom(\Ga,G)$, leaving $R(\Ga,G,H)$ invariant. We now
define the following two spaces:
\begin{description}
\item[the deformation space] $\sT(\Ga,G,H):=R(\Ga,G,H)/G$; and
\item[the moduli space] $\sM(\Ga,G,H):=\Aut(\Ga)\bs R(\Ga,G,H)/G$.
\end{description}
For example, if $(G,H)=(\PSL(2,\R),\PSO(2))$ and $\Ga=\pi_1(M_g)$ for
$g\ge2$ then $\sT(\Ga,G,H)$ is the Teichm\"uller space of $M_g$, and
$\sM(\Ga,G,H)$ is nothing other than the moduli space of complex
structures on $M_g$.

We formalize local rigidity of a discontinuous group as saying that it
corresponds to an ``isolated point'' of the deformation space
$\sT(\Ga,G,H)$:

\begin{dfn}[Local rigidity in a non-Riemannian homogeneous space 1993
\cite{20}]\rm Let $u\in R(\Ga,G,H)$. We say that the discontinuous group
$u(\Ga)$ for the homo\-geneous space $G/H$ determined by $u$ is {\em
locally rigid} as a discontinuous group of $G/H$ if the single point
$[u]$ is an open set of the quotient space $\Hom(\Ga,G)/G$; this means
that any point sufficiently close to $u$ is conjugate to $u$ under an
inner automorphism of $G$. If $u$ is not locally rigid, we say that $u$
{\em admits continuous deformations}. \end{dfn}

When $H$ is compact, this terminology coincides with the original notion
(see Weil \cite{47}).

In higher dimensions, let us compare whether the local rigidity theorem
holds in the cases that $H$ is compact or noncompact. Let $G'$ be a
non\-compact simple linear Lie group, and $K'$ its maximal compact
subgroup. We can use the vanishing and non\-vanishing theorems for
cohomology of Lie algebras together with the criterion for a properly
discontinuous action discussed above (Theorem~14) and so on, to prove
the following theorem.

\begin{thm} {\bf (Local rigidity theorem -- the Riemannian case: 
\\
Selberg
and Weil 1964 \cite{47})}\enspace
Let $(G,H):=(G',K')$. Then the following two conditions on $G'$ are
equi\-valent:
\begin{enumerate}
\renewcommand{\labelenumi}{(\roman{enumi})}
\item there exists a uniform lattice $\iota\colon\Ga\to G'$ such that
$\iota\in R(\Ga,G,H)$ admits continuous deformations.
\item $G'$ is locally isomorphic to $\SL(2,\R)$.
\end{enumerate}
\end{thm}

\begin{thm}[Local rigidity theorem -- the non-Riemannian case 1998 \cite{24}]
Let $(G,H):=(G'\times G',\diag(G'))$. Then the following two conditions on $G'$ are
equivalent:
\begin{enumerate}
\renewcommand{\labelenumi}{(\roman{enumi})}
\item there exists a uniform lattice $\iota\colon\Ga\to G'$ such that
$\iota\times1\in R(\Ga,G,H)$ admits continuous deformations.
\item $G'$ is locally isomorphic to $\SO(n,1)$ or $\SU(n,1)$.
\end{enumerate}
\end{thm}

From a different point of view, for the group manifold $G'$ and its
uniform lattice $\Ga$, Theorem~27 only treats the rigidity of left
actions, whereas Theorem~28 considers the rigidity of both left and
right actions. (Note that $G/H\cong G'$ in the latter case.)

The deformation spaces have also been studied in the following cases:

\begin{enumerate}
\renewcommand{\labelenumi}{(\arabic{enumi})}
\item The Poincar\'e disk $G/H=\SL(2,\R)/\SO(2)$.
\item $G/H=G'\times G'/\diag G'$ for $G'=\SL(2,\R)$ (Goldman 1985 \cite{11},
Salein 1999 \cite{45}).
\item $G/H=G'\times G'/\diag G'$ for $G'=\SL(2,\C)$ (Ghys 1995 \cite{13}).
\end{enumerate}

In these cases the deformation space $\sT(\Ga,G,H)$ corresponds
respectively to:
\begin{enumerate}
\renewcommand{\labelenumi}{(\arabic{enumi})}
\item The deformations of complex structures on a Riemann surface $M_g$
of genus $g\ge2$.
\item The deformations of Lorentz structures on a 3-dimensional
manifold.
\item The deformation of complex structures on a 3-dimensional complex
manifolds.
\end{enumerate}
Moreover, (2) and (3) correspond to the cases $n=1,2,3$ of Theorem~28.
Indeed, this follows because we have the local isomorphisms of Lie
groups
\[
G' = \SL(2,\R) \approx \SO(2,1) \approx \SU(1,1)
\quad\hbox{and}\quad
G' = \SL(2,\C) \approx \SO(3,1).
\]

In Theorem~28, as $n$ increases, one sees that one can construct
irreducible pseudo-Riemannian symmetric spaces of arbitrarily high
dimension, together with a uniform lattice $\Ga$ for which the local
rigidity theorem does not hold. In \cite{24}, for general $n$, we
obtained quantitative estimates for deformations of this type of uniform
lattice $\Ga$ (that is, the range within which discontinuity is
preserved) using the diameter of locally Riemannian symmetric spaces
$\Ga\bs G'/K'$. Now the proposition
\begin{quote}
a ``small'' deformation of a discrete subgroup preserves the
discontinuity of the action
\end{quote}
is false for general Lie groups (see \cite{29}).%
However, in the case of semisimple Lie groups discussed above this
proposition holds, and in particular this settles Goldman's conjecture
\cite{11} positively. The key to the proof is the criterion for a
properly discontinuous action (Theorem~14).

\bigskip

\noindent
Received July 16th 2004
\bigskip

\noindent
KOBAYASHI Toshiyuki, Kyoto University,\\
Research Institute for Mathematical Sciences,\\
e-mail \verb!toshi@kurims.kyoto-u.ac.jp!

\bigskip

\noindent
Translated by M. Reid\\
\verb!miles@maths.warwick.ac.uk!

\end{document}